\documentclass[a4paper,12pt]{article}
\usepackage{amssymb,amsmath}
\usepackage{latexsym}
\usepackage{graphicx}
\newtheorem{theorem}{Theorem}

\newtheorem{conjecture}{Conjecture}

\newtheorem{question}{Question}

\def\bc{\begin{center}}\def\ec{\end{center}}
\def\beq{\begin{equation}}\def\eeq{\end{equation}}
\def\beqn{\begin{eqnarray}}\def\eeqn{\end{eqnarray}}
\def\bd{\begin{eqnarray*}}\def\ed{\end{eqnarray*}}
\def\pont{\hspace{-6pt}{\bf.\ \ }}
\def\:{\!:}
\def\qed{\ifhmode\unskip\nobreak\fi\quad\ifmmode\Box\else$\Box$\fi}

\title{Rainbow matchings in bipartite multigraphs}
\author{J\'anos Bar\'at \\
\small  MTA-ELTE Geometric and Algebraic Combinatorics Research Group\\[-0.8ex]
\small \texttt{barat@cs.elte.hu} \and
Andr\'as Gy\'arf\'as\thanks{Research was supported in part by OTKA K104343.}\\[-0.8ex]
\small Alfr\'ed R\'enyi Institute of Mathematics\\[-0.8ex]
\small Hungarian Academy of Sciences\\[-0.8ex]
\small Budapest, P.O. Box 127\\[-0.8ex]
\small Budapest, Hungary, H-1364\\[-0.8ex] \small
\texttt{gyarfas.andras@renyi.mta.hu} \and G\'{a}bor N. S\'ark\"ozy\thanks{Research was supported in part by OTKA K104343.}\\[-0.8ex]
\small Alfr\'ed R\'enyi Institute of Mathematics\\[-0.8ex]
\small Hungarian Academy of Sciences\\[-0.8ex]
\small Budapest, P.O. Box 127\\[-0.8ex]
\small Budapest, Hungary, H-1364\\[-0.8ex]
\small \texttt{sarkozy.gabor@renyi.mta.hu}\\
and\\
\small Computer Science Department\\[-0.8ex]
\small Worcester Polytechnic Institute\\[-0.8ex]
\small Worcester, MA, USA 01609\\[-0.8ex]
\small \texttt{gsarkozy@cs.wpi.edu}}

\date{\today}

\begin{document}

\maketitle

\begin{abstract}

Suppose that $k$ is a non-negative integer and a bipartite
multigraph $G$ is the union of $$N=\left\lfloor
\frac{k+2}{k+1}n\right\rfloor -(k+1)$$ matchings $M_1,\dots,M_N$,
each of size $n$. We show that $G$ has a rainbow matching of size
$n-k$, i.e. a matching of size $n-k$ with all edges coming from
different $M_i$'s. Several choices of the parameter $k$ relate to known
results and conjectures.
\end{abstract}


Suppose that a multigraph $G$ is given with a proper $N$-edge
coloring, i.e. the edge set of $G$ is the union of $N$ matchings
$M_1,\dots,M_N$. A {\em rainbow matching} is a matching whose edges
are from different $M_i$'s.

A well-known conjecture of Ryser \cite{RY} states that for odd $n$
every $1$-factorization of $K_{n,n}$ has a rainbow matching of size
$n$. The companion conjecture, attributed to Brualdi \cite{BR} and
Stein \cite{S} states that for every $n$, every $1$-factorization of
$K_{n,n}$ has a rainbow matching of size at least $n-1$. These
conjectures are known to be true in an asymptotic sense, i.e. every
$1$-factorization of $K_{n,n}$ has a rainbow matching containing
$n-o(n)$ edges. For the $o(n)$ term, Woolbright \cite{WO} and
independently Brouwer et al. \cite{BVW} proved $\sqrt{n}$. Shor
\cite{SH} improved this to $5.518(\log{n})^2$, an error was
corrected in \cite{HS}.




There are several results for the case when $K_{n,n}$ is replaced by an arbitrary
bipartite multigraph. The following conjecture of Aharoni et al.
\cite{ACH} strengthens the Brualdi-Stein conjecture.

\begin{conjecture}\label{aharoni}\pont
If a bipartite multigraph $G$ is the union of $n$ matchings of size
$n$, then $G$ contains a rainbow matching of size $n-1$.
\end{conjecture}

As a relaxation,
Kotlar and Ziv \cite{KZ} noticed that the union of $n$ matchings of size
$\frac{3}{2}n$ contains a rainbow matching of size $n-1$.
Conjecture~\ref{aharoni} would follow from another one posed by Aharoni and Berger:

\begin{conjecture}\label{unpublished}\pont
If a bipartite multigraph $G$ is the union of $n$ matchings of size
$n+1$, then $G$ contains a rainbow matching of size $n$.
\end{conjecture}

Recently, there has been gradual progress on this question.
Aharoni et al. proved that matchings of size $\frac{7}{4}n$ suffice \cite{ACH}.
Kotlar and Ziv \cite{KZ} improved it to $\frac{5}{3}n$ and Clemens and Ehrenm\"uller to $(\frac{3}{2}+\varepsilon)n$.

One needs a lot more matchings of size $n$ to guarantee a rainbow
matching of size $n$.  Aharoni and Berger \cite{AB} and (in a
slightly weaker form) Drisko \cite{D} proved the following.

\begin{theorem} \label{drisko}\pont
If a bipartite multigraph $G$ is the union of $2n-1$ matchings of
size $n$, then $G$ contains a rainbow matching of size $n$.
\end{theorem}

The (unique) factorization of a cycle on $2n$ vertices with edges of
multiplicity $n-1$ shows that in the statement $2n-1$ cannot be
replaced by $2n-2$ (see \cite{D}). We merge Conjecture \ref{aharoni}
and Theorem \ref{drisko} into a unified context and ask the
following. (We note that this question was also raised independently
in \cite{CE}.)

\begin{question}\pont\label{question}
For integers $0\leq k < n$, what is the smallest $N=N(n,k)$ such
that any bipartite multigraph $G$ that is the union of $N$ matchings
of size $n$, contains a rainbow matching of size $n-k$?
\end{question}

Conjecture \ref{aharoni} claims that $N(n,1)=n$ and Theorem
\ref{drisko} states that $N(n,0)=2n-1$.
 In this note we give the following
upper bound on $N(n,k)$.

\begin{theorem}\label{bip}\pont For $0\leq k < n$,
$N(n,k)\le \left\lfloor \frac{k+2}{k+1}n\right\rfloor -(k+1)$.
\end{theorem}

In the range $\lfloor n/2\rfloor\le k < n$ Theorem \ref{bip} gives
$N(n,k)\le n-k$ which is obviously best possible, therefore
$N(n,k)=n-k$. When $k=0$ it gives $N(n,0)\le 2n-1$, the bound of
Theorem \ref{drisko}, so this is best possible as well. The case
$k=1$ gives a result towards Conjecture~\ref{aharoni}: if a
bipartite multigraph is the union of $\lfloor\frac{3}{2}n\rfloor-2$
matchings of size $n$, then there is a rainbow matching of size
$n-1$. As far as we know this is the best result in this direction.
If $N=\lfloor(1+\epsilon)n\rfloor$ for some $\epsilon>0$, we get a
partial rainbow matching of size $n-c$ where $c$ is a constant
depending on $\epsilon$ ($c=\lfloor 1/\epsilon\rfloor$), this goes
beyond the best error term known for Ryser's conjecture (\cite{HS}),
but the price is the increment in the number of colors. Also, when
$k=\lfloor \sqrt{n}\rfloor$, Theorem \ref{bip} extends (from
factorizations of $K_{n,n}$ to colorings of bipartite multigraphs)
Woolbright's result \cite{WO}, namely that a factorization of
$K_{n,n}$ contains a rainbow matching of size at least $n-\sqrt{n}$.

\bigskip

\noindent {\bf Proof of Theorem \ref{bip}. \rm}  We use Woolbright's argument \cite{WO}.
Set \mbox{$N=\left\lfloor \frac{k+2}{k+1}n\right\rfloor -(k+1)$}.
Let the edge set of a bipartite multigraph
$G=[A,B]$ be the union of matchings $M_1,\dots,M_N$ each of size
$n$ and let $R_1$ be a maximum rainbow matching of $G$
with $t$ edges. Suppose to the contrary that $t\le n-k-1$.

We assume the edges of $M_1,\dots,M_{N-t}$ are not
used in $R_1$. For any subset $S\subset B$, define
$$f(S)=\{v\in A: (v,w)\in R_1 \mbox{ for some } w\in S\}.$$

Set  $B_0=B\setminus V(R_1),A_0=A\setminus V(R_1)$.  For every $j\in
\{1,\dots,N-t\}$ a matching $F_j\subset M_j$ of size $j(n-t)$ will
be defined with the following property.
\begin{itemize}
\item Property 1: $V(F_j)\cap B_0=\emptyset.$
\end{itemize}

Let $F_1\subset M_1$ be a matching of size $n-t$ such that $V(F_1)\cap A \subseteq A_0$, since $|M_1|-|R_1|=n-t$, such $F_1$ exists.
Set $B_1=V(F_1)\cap B$.
Since $R_1$ is a maximum rainbow matching, $V(F_1)\cap B_0=\emptyset$, so Property 1 holds and $|F_1|=1\times (n-t)$.
Set $A_1=f(B_1)$.

Suppose that for some $i\ge 1$ the matchings $F_i,R_i$ and the pairwise disjoint $(n-t)$-element sets
$A_1,\dots,A_i,B_1, \ldots,B_i$ have already been defined, where $|F_i|=i(n-t)$.
Define the rainbow matching $R_{i+1}$ by removing from $R_i$ the edges that go from $B_i$ to $A_i$.

To define $F_{i+1}\subset M_{i+1}$, take $(i+1)(n-t)$ edges of $M_{i+1}$
incident to $A\setminus V(R_{i+1})$.
There exist sufficiently many edges in $M_{i+1}$ since
$$|M_{i+1}|-|R_{i+1}|=n-(t-\sum_{j=1}^i |B_j|)=(i+1)(n-t).$$

We show that Property 1 is maintained.
Suppose to the contrary that we find $(a_0,b_0)\in
F_{i+1}$, $a_0\in A_{j}$ for some $1\leq j \leq i$, $b_0\in B_0$
(clearly $j\not=0$). Then $b_1=f^{-1}(a_0)\in B_j$, and there
exists an $a_1$ such that $(a_1,b_1)\in F_j$  and this generates an
alternating path
$$Q=(b_0,a_0),(a_0,f^{-1}(a_0)),(f^{-1}(a_0),a_1),(a_1,f^{-1}(a_1)),(f^{-1}(a_1),a_2),\dots$$
ending in $A_0$ allowing us to replace all edges of $R_1\cap E(Q)$
by edges in different $F_j$s ($j\le i+1$) contradicting the choice
of $t$. Note that $Q$ is a simple path, since with some
$j>j_1>\dots>j_k>0$, its edges go between the disjoint sets
$$(B_0,A_j),(A_j,B_j),(B_j,A_{j_1}),(A_{j_1},B_{j_1}),(B_{j_1},A_{j_2}),\dots,(A_{j_k},B_{j_k}),(B_{j_k},A_0).$$

Now $F_{i+1}$ is defined and by Property 1
$$|V(F_{i+1})\cap (B\setminus (\cup_{k=0}^i B_k))|\ge n-t,$$
therefore we can define $B_{i+1}$ as an $(n-t)$-element subset of
$V(F_{i+1})\cap (B\setminus (\cup_{k=0}^i B_k))$. Finally, set
$A_{i+1}=f(B_{i+1})$.

Since $V(F_{N-t})\cap B \subseteq B\setminus B_0$, we get
$$(N-t)(n-t)\le t.$$
Dividing by $n-t$ (using $t\le n-k-1<n$) this can be rewritten as
$$N-t\leq \frac{t}{n-t}=\frac{n-n+t}{n-t} = \frac{n}{n-t}-1$$
or
$$N\leq \frac{n}{n-t}+t-1.$$
Using this, the definition of $N$ and $t\le n-k-1$, we get
$$\left\lfloor \frac{k+2}{k+1}n\right\rfloor -(k+1)=N\le \frac{n}{n-t}+t-1\le \frac{n}{k+1}+n-k-1-1
$$
and this leads to
$$\left\lfloor \frac{n}{k+1}\right\rfloor \le \frac{n}{k+1}-1,$$
a contradiction, finishing the proof. \qed

\noindent {\bf Remark.} A natural variant of Question \ref{question}
is to allow arbitrary multigraphs (instead of bipartite ones).
Denote the corresponding function by $N'(n,k)$. For $k=0$ we have an
example showing $N'(n,0)>2n-1$ and recently Aharoni informed us
\cite{A} that they proved $N'(n,0)\le 3n-2$. Indeed, our example is
the following. Let the vertices be denoted as $1,2,\dots,4k$, where
$2n=4k$. Let $M_1=\dots=M_{n-1}=\{12, 34, \dots, (2n-1)2n\}$,
$M_{n}=\dots=M_{2n-2}=\{23, 45, \dots, (2n)1\}$ and
$M_{2n-1}=\{13,24,57,68,\dots,(2n-3)(2n-1), (2n-2)2n\}$. As it was
remarked before, there is no full rainbow matching without using an
edge of $M_{2n-1}$. We may assume that we use the edge $24$. Now any
edge of $M_i$ that covers the vertex $3$, where $1\le i\le 2n-2$,
uses either vertex $2$ or $4$. Therefore, there is no full rainbow
matching.

\end{document}